\documentclass[a4paper,10pt]{article}

\usepackage{amsmath}
\usepackage{ntheorem}
\usepackage{eucal}
\usepackage{amscd}
\usepackage{ifthen}
\usepackage{amssymb}
\usepackage[english,french]{babel}
\usepackage[T1]{fontenc}
\usepackage[latin1]{inputenc}
\usepackage[all]{xy}

\newenvironment{demo}[1][]{\ifthenelse{\equal{#1}{}}{\noindent \textbf{D\'emonstration :\\ \indent}}{\noindent \textbf{D\'emonstration #1 :\\ \indent}}}{$\square$\\}

\theoremstyle{break}
\newtheorem{theo}{Théorème}[section]
\newtheorem{cor}{Corollaire}[section]

\newtheorem{lem}{Lemme}[section]
\newtheorem{prop}{Proposition}[section]
\newtheorem{defi}{Définition}[section]
\theorembodyfont{\slshape}
\newtheorem{rem}{Remarque}[section]
\theorembodyfont{\upshape}

\newcommand{\RR}{\mathbb R}

\newcommand{\CC}{\mathbb C}
\newcommand{\QQ}{\mathbb Q}
\newcommand{\ZZ}{\mathbb Z}

\newcommand{\To}{\longrightarrow}

\newcommand{\set}[1]{\left\{#1\right\}}

\newcommand{\hg}[2]{\pi_{#1}(#2)}
\newcommand{\dimm}[1]{\mathrm{dim}(#1)}
\newcommand{\alb}[1]{\mathrm{Alb}(#1)}
\newcommand{\albs}[2]{\mathrm{Alb}_{#1}(#2)}

\newcommand{\wtx}{\widetilde{X}}

\newcommand{\tor}[1]{\mathrm{Tor}\left(#1\right)}

\newcommand{\shm}{\textsc{shm}~}

\title{Convexité holomorphe du revêtement de Mal\v{c}ev d'après Sandrine Leroy}
\author{Benoît \textsc{Claudon}}
\date{}

\begin{document}

\maketitle

\begin{abstract}
Dans cette note, nous présentons l'approche de \cite{L98} concernant le cas nilpotent de la conjecture de Shafarevich. Le cas projectif est traité dans \cite{K97} mais la démonstration de S. Leroy procède de manière indépendante et s'applique aussi bien au cas kählérien. Elle est entièrement fondée sur l'existence des variétés d'Albanese supérieures (construites par R. Hain).
\end{abstract}
\vspace*{0.5cm}

\section*{Introduction}

La conjecture de Shafarevich \cite{Sh74} prédit que le revêtement universel $\wtx$ d'une variété projective $X$ est holomorphiquement convexe (\emph{i.e.} admet une application holomorphe propre sur un espace de Stein). Dans toute sa généralité, cette conjecture est actuellement hors de portée mais, en imposant certaines conditions algébriques sur le groupe fondamental de $X$, il devient possible d'obtenir la convexité holomorphe de $\wtx$. Ainsi, si le groupe fondamental est presque-abélien\footnote{un groupe est dit presque-abélien (resp. presque-nilpotent) si il admet un sous-groupe d'indice fini abélien (resp. nilpotent).}, on montre aisément que $\wtx$ admet une application propre sur $\CC^n$ (revêtement universel de $\alb{X'}$ avec $X'\To X$ revêtement étale fini de $X$) et $\wtx$ vérifie donc la conjecture de Shafarevich.

Dans la catégorie des groupes de type fini, les groupes nilpotents succèdent immédiatement aux groupes abéliens en ce qui concerne la complexité algébrique. Nous montrons ici que la conjecture de Shafarevich est encore vérifiée pour cette catégorie de groupes (fondamentaux). Plus précisément, si
\begin{equation*}
\hg{1}{X}\To \hg{1}{X}_{nilp}
\end{equation*}
désigne la complétion nilpotente (modulo torsion) de $\hg{1}{X}$ (voir la section \ref{rappels} pour les notions utilisées), on notera $X^{nilp}$ le revêtement dont le groupe de Galois est $\hg{1}{X}_{nilp}$ ; $X^{nilp}$ est aussi appelé revêtement de Mal\v{c}ev.

\begin{theo}[\cite{K97},\cite{L98}]\label{shaf nilpotent}
Le revêtement de Mal\v{c}ev d'une variété kählérienne compacte $X$ est holomorphiquement convexe.
\end{theo}

\begin{cor}\label{cor shaf}
Si $\hg{1}{X}$ est presque-nilpotent, $\wtx$ est holomorphiquement convexe.
\end{cor}

Avant d'esquisser les grandes lignes de la démonstration de \cite{L98}, rappelons la façon dont procède \cite{K97} dans la situation d'un groupe $\pi_1(X)$ nilpotent sans torsion. Un tel groupe admettant des quotients abéliens libres non triviaux, l'application d'Albanese de $X$ a une image de dimension strictement positive ; notons $\alpha:X\To S$ la factorisation de Stein de cette application. Remarquons tout de suite que si $S$ est lisse, la démonstration du théorème \ref{shaf nilpotent} est élémentaire et est reproduite au paragraphe \ref{cas part}. Toute la difficulté réside donc dans les éventuelles singularités de $S$, cet obstacle étant surmonté grâce à l'utilisation de Structures de Hodge Mixtes\footnote{\textsc{shm}  dans toute la suite.} (voir le paragraphe \ref{shm} pour les notions considérées).

Les fibres de $\alpha$ sont les sous-variétés $Z$ de $X$ vérifiant :
\begin{equation}\label{fibres alb}
\mathrm{Im}\left(H_1(Z,\QQ)\To H_1(X,\QQ)\right)=0.
\end{equation}
L'existence de \shm  sur les algèbres de Lie $L_\CC(\pi_1(X))$ et $L_\CC(\pi_1(Z))$ (théorème \ref{shm}) et leur caractère fonctoriel, combinés avec le caractère \emph{strict} des morphismes de \shm  (voir \cite{Del71}), permet de passer de l'homologie à l'homotopie :
\begin{equation}\label{homologie-homotopie}
\mathrm{Im}\left(H_1(Z,\QQ)\To H_1(X,\QQ)\right)=0\Longleftrightarrow \mathrm{Im}\left(\pi_1(Z)\To\pi_1(X)\right)\textrm{ est fini.}
\end{equation}
L'application $\alpha$ induit alors une application
$$\tilde{\alpha}:\wtx\To \widetilde{S}$$
entre les revêtement universels qui est \emph{propre} en vertu de (\ref{homologie-homotopie}). D'autre part, $S$ étant finie sur un tore, son revêtement universel $\widetilde{S}$ est Stein et $\wtx$ est bien holomorphiquement convexe. Remarquons enfin que l'hypothèse projective est essentielle dans la démonstration de \cite{K97}, la démonstration de (\ref{homologie-homotopie}) se faisant par réduction aux cas des courbes (grâce au théorème de Lefschetz).\
\vspace*{0.3cm}

La démonstration présentée ici s'appuie de manière essentielle sur l'existence des variétés d'Albanese supérieures. Ces variétés, construites par R. Hain (voir \cite{H85}), jouent le même rôle vis-à-vis des quotients nilpotents (sans torsion) de $\hg{1}{X}$ que la variété d'Albanese usuelle $\alb{X}$ pour l'abélianisé du groupe fondamental (qui n'est autre que le premier quotient de la tour nilpotente). Ainsi, si $(\albs{s}{X})_{s\geq 1}$ désigne la suite des variétés d'Albanese supérieures, le groupe fondamental
\begin{displaymath}
\hg{1}{\albs{s}{X}}
\end{displaymath}
est en particulier canoniquement isomorphe au $s^{\textrm{ième}}$ quotient sans torsion de $\hg{1}{X}$. Remarquons tout de suite que ces variétés sont essentiellement des quotients de groupes de Lie (simplement connexes) complexes nilpotents et qu'elles ne sont en général pas compactes, ni même kählériennes. Ces variétés s'organisent en une tour de fibrations principales :
$$\cdots\To\albs{s+1}{X}\To \albs{s}{X}\To\cdots\To \albs{1}{X}=\alb{X}$$
de fibres des quotients de groupes de Lie complexes abéliens. Rappelons de plus que la variété $X$ est elle-même munie de morphismes
\begin{displaymath}
\alpha_s:X\To\albs{s}{X}
\end{displaymath}
qui induisent les isomorphismes évoqués ci-dessus.

La démonstration consiste alors à remarquer que la suite des images
$$(\alpha_s(X)\subset\albs{s}{X})_{s\geq1}$$
stationne pour $s$ assez grand (ce qui n'est bien entendu pas nécessairement le cas de la tour nilpotente). Le revêtement de Mal\v{c}ev $X^{nilp}$ est alors essentiellement obtenu en tirant en arrière par l'application $\alpha_s$ le revêtement universel de $\albs{s}{X}$ (pour $s$ assez grand). Le caractère Stein de ce dernier permet alors de conclure quant à la convexité holomorphe de $X^{nilp}$.\\

Dans un premier temps, nous donnerons quelques éléments de rappels sur la suite centrale d'un groupe de type fini ainsi que sur la complétion de Mal\v{c}ev d'un groupe (de type fini) nilpotent sans torsion. L'existence d'une \shm  sur la complétion de Mal\v{c}ev de $\hg{1}{X}$ (pour $X$ variété kählérienne compacte) permet alors d'expliciter la construction des variétés d'Albanese supérieures et d'en obtenir les propriétés essentielles. Dans la dernière partie, nous montrerons comment rendre effectives les grandes lignes décrites ci-dessus.

\begin{rem}
Après avoir étudié le cas nilpotent de la conjecture de Shafarevich, il semblerait naturel de s'intéresser à la situation d'un groupe fondamental résoluble. Or, il se trouve que la catégorie des groupes kählériens possèdent des propriétés très particulières ; en particulier, T. Delzant montre dans \cite{D07} qu'un groupe kählérien résoluble est automatiquement presque-nilpotent. Le cas résoluble de la conjecture de Shafarevich est donc compris dans le corollaire \ref{cor shaf}.
\end{rem}

\begin{rem}
Concernant la conjeture de Shafarevich, les résultats les plus significatifs à l'heure actuelle concernent les revêtement linéaires, comme c'est le cas dans \cite{E04}. Les groupes nilpotents (de type fini) étant eux-mêmes linéaires (voir \cite[chap. 5]{S83}), les résultats de \cite{E04} impliquent donc le théorème \ref{shaf nilpotent}. Il nous a cependant semblé intéressant d'exposer la démonstration ci-dessous pour son caractère original.
\end{rem}

\section{Rappels et préliminaires}\label{rappels}

\subsection{Suite centrale descendante}

Soit $G$ un groupe de type fini. On associe à $G$ sa suite centrale decendante définie par :
\begin{displaymath}
G_1=G\quad\textrm{et}\quad\forall\,i\geq1,\,G_{i+1}=[G_i,G]
\end{displaymath}
Un tel groupe $G$ est dit nilpotent de classe au plus $s$ si $G_{s+1}=\set{1}$. Les quotients successifs associés à la suite centrale s'organise en une tour
\begin{displaymath}
\cdots G/G_{i+1}\To G/G_i\To \cdots G/G_3\To G/G_2=G_{ab}\To 1,
\end{displaymath}
les extensions
\begin{displaymath}
1\To G_i/G_{i+1}\To G/G_{i+1}\To G/G_i\To 1
\end{displaymath}
étant centrales par construction, avec $G/G_{i+1}$ nilpotent de classe au plus $i$.

\begin{prop}[lem. 1.2, p. 470 \cite{P77}]
On note
\begin{displaymath}
\tor{G}=\set{g\in G\vert\,\exists\, n\geq1,\, g^n=1}
\end{displaymath}
la torsion de $G$. Si $G$ est nilpotent, $\tor{G}$ est un sous-groupe caractéristique de $G$ et le groupe $G_{st}=G/\tor{G}$ est sans torsion. Si de plus $G$ est de type fini, $\tor{G}$ est fini.
\end{prop}
Ceci s'applique en particulier aux quotients $G/G_i$ ; on note alors
\begin{displaymath}
G'_i=\mathrm{Ker}\left(G\To \left(G/G_i\right)_{st}\right)=\set{g\in G\vert\,\exists\, n\geq1,\, g^n\in G_i}.
\end{displaymath}

\begin{defi}
Soit $G$ un groupe de type fini. On note
$$G_\infty=\bigcap_{i\geq1}G'_i\textrm{ et }G_{nilp}=G/G_\infty.$$
$G_{nilp}$ s'appelle la complétion nilpotente modulo torsion ; c'est le plus gros quotient de $G$ obtenu comme limite projective de groupes nilpotents sans torsion.
\end{defi}

\subsection{Complétion de Mal\v{c}ev}

La catégorie des groupes nilpotents de type fini possède des propriétés de finitude tout à fait remarquables. Ainsi, si $G$ est un groupe de type fini, les quotients $G_i/G_{i+1}$ restent de type fini \cite[cor. 7, p. 13]{S83} et, si $G$ est de plus nilpotent, il est obtenu par une suite finie d'extensions centrales par des groupes de type fini.

La suite centrale descendante possède une autre propriété intéressante, à savoir :
\begin{equation}\label{prop lie}
\forall\,i,j\geq1,\,[G_i,G_j]\subset G_{i+j}.
\end{equation}
Si $K$ désigne un corps de caractéristique 0 et si $G$ est nilpotent (de classe $s$) de type fini, on considère
$$L_K(G)=\bigoplus_{i=1}^s\left(G_i/G_{i+1}\right)\otimes K$$
qui est canoniquement munie d'une structure d'algèbre de Lie (de dimension finie) sur $K$ : d'après \ref{prop lie}, le crochet $[x,y]=xyx^{-1}y^{-1}$ dans $G$ induit un effet un crochet de Lie sur $L_K(G)$. Celle-ci hérite immédiatement des propriétés de $G$ : c'est une algèbre de Lie nilpotente. La formule de Campbell-Haussdorf se réduit alors en une identité polynômiale et permet de munir l'espace vectoriel sous-jacent à $L_K(G)$ d'une structure de groupe de Lie dont l'algèbre de Lie est $L_K(G)$ (voir le chapitre 5 de \cite{BK81}). Le groupe ainsi obtenu est une sorte de complétion de $G$.
\begin{theo}[Mal\v{c}ev, \cite{M49}]\label{complétion Malcev}
Soit $G$ un groupe nilpotent de type fini et $K$ un corps de caractéristique 0. Il existe alors un unique groupe de Lie $\hat{G}_K$ nilpotent défini sur $K$ et d'algèbre de Lie $L_K(G)$, ainsi qu'un morphisme
$$i:G\To \hat{G}_K$$
tel que :
\begin{enumerate}
\item $\textrm{Ker}(i)=\textrm{Tor}(G)$ et
\item $\textrm{Im}(i)$ est cocompacte dans $\hat{G}_\RR$.
Si $G$ est sans-torsion, on peut donc canoniquement réaliser $G$ comme réseau d'un groupe de Lie nilpotent (simplement connexe).
\end{enumerate}
\end{theo}
Le groupe $\hat{G}_K$ s'appelle la complétion de Mal\v{c}ev de $G$. Dans la suite, on va s'intéresser au cas où $G$ est un des quotients nilpotents sans-torsion de $\hg{1}{X}$ avec $X$ kählérienne compacte.

\section{Variétés d'Albanese supérieures}

\subsection{Structures de Hodge mixtes}\label{shm}

Nous allons brièvement exposer ici la construction des variétés d'Albanese supérieures supérieures (dûe à R. Hain) associées à une variété kählérienne compacte $X$. Celle-ci repose sur l'existence d'une \shm sur la complétion de Mal\v{c}ev du groupe fondamental de $X$. Afin d'alléger le formalisme introduit ci-dessus, nous adopterons les notations suivantes pour $X$ est une variété kählérienne compacte et $s\geq 1$ un entier :
\begin{description}
\item[-] $\mathcal{G}_s(\CC)=L_\CC\left(\hg{1}{X}/\hg{1}{X}'_{s+1}\right)$ et
\item[-] $G_s(\CC)$ la complétion de Malcev de $\hg{1}{X}/\hg{1}{X}'_{s+1}$ définie sur $\CC$.
\end{description}

\noindent Rappelons tout d'abord la notion de \shm.
\begin{defi}
Une \shm de poids $k$ sur un $\QQ$-espace vectoriel de dimension finie $H_\QQ$ est la donnée d'une filtration croissante $W^{\bullet}$ de $H_\QQ$ (dite filtration par le poids) et d'une filtration décroissante $F^{\bullet}$ de $H_\CC=H_\QQ\otimes\CC$ (dite filtration de Hodge) compatibles dans le sens suivant : la filtration $F^{\bullet}$ induit sur chaque gradué
$$\mathrm{Gr}^W_n(H)=W^{n}(H_\QQ)/W^{n-1}(H_\QQ)$$
une structure de Hodge (rationnelle) pure de poids $n+k$.
\end{defi}

\noindent L'un des intérêts majeurs des \shm réside dans le résultat suivant.
\begin{theo}[Deligne, \cite{Del71}]
Les morphismes de \shm sont \emph{stricts} pour les filtrations $W^{\bullet}$ et $F^{\bullet}$ : si
$$\phi:(H_\QQ,W,F)\To (V_\QQ,U,G)$$
est un morphisme de \shm, on a alors
\begin{eqnarray*}\mathrm{Im}(\phi)\cap U^{\bullet}(V_\QQ)&=&\phi(W^{\bullet}(H_\QQ))\\
\mathrm{Im}(\phi)\cap G^{\bullet}(V_\CC)&=&\phi(F^{\bullet}(H_\CC))
\end{eqnarray*}
\end{theo}

\noindent Les travaux de Hain montrent que l'algèbre de Lie de Mal\v{c}ev du groupe fondamental d'une variété kählérienne compacte est naturellement munie d'une \shm.
\begin{theo}[Hain, \cite{H87}]
Soit $X$ une variété kählérienne compacte et $s\geq1$. Il existe une \shm sur $\mathcal{G}_s(\CC)$ (fonctorielle en $X$) dont la filtration par le poids $W^\bullet$ est donnée par la suite centrale descendante de l'algèbre de Lie $\mathcal{G}_s(\CC)$. De plus, pour cette structure, le crochet de Lie est un morphisme de structure de Hodge mixte.
\end{theo}

La description complète de cette \shm est malheureusement trop élaborée pour être exposée en quelques lignes ; nous renvoyons à \cite{PS08} pour une exposition raisonnable des constructions de Hain.

Notons $F^\bullet$ la filtration de Hodge obtenue sur $\mathcal{G}_s(\CC)$ et considérons l'application exponentielle
$$\exp:\mathcal{G}_s(\CC)\To G_s(\CC).$$
Le groupe de Lie $G_s(\CC)$ étant nilpotent et simplement connexe, cette application est un difféomorphisme et, par abus de notation, nous noterons
$$F^0G_s(\CC)=\exp(F^0(\mathcal{G}_s(\CC)))$$
le sous-groupe fermé associée à la sous-algèbre $F^0(\mathcal{G}_s(\CC))$.

\subsection{Construction des variétés d'Albanese supérieures}

Dans \cite{H85}, R. Hain montre le résultat suivant :
\begin{theo}[Hain, \cite{H85}]\label{albanese sup}
Soit $X$ une variété kählérienne compacte. Il existe alors des variétés lisses $\left(\albs{s}{X}\right)_{s\geq1}$ et des morphismes (holomorphes)
$$\alpha_s:X\To\albs{s}{X}\textrm{ et }\pi_s:\albs{s+1}{X}\To\albs{s}{X}$$
rendant le diagramme suivant commutatif :
$$\xymatrix{
&&X \ar[ld]_{\alpha_{s+1}}\ar[d]^{\alpha_s}\ar[rrd]^{\alpha_1}&&\\
\ar@{-->}[r]&\albs{s+1}{X}\ar[r]^{\pi_s}&\albs{s}{X}\ar@{-->}[rr]&&\albs{1}{X}
}$$
Les propriétés suivantes sont également satisfaites :
\begin{enumerate}
\item $\albs{1}{X}$ coincide avec la variété d'Albanese usuelle (ainsi que $\alpha_1$),
\item les morphismes $\alpha_s$ induisent des isomorphismes
$$\alpha_s:\hg{1}{X}/\hg{1}{X}'_s\stackrel{\sim}{\To}\hg{1}{\albs{s}{X}}$$
\item et les projections $\pi_s:\albs{s+1}{X}\To\albs{s}{X}$ sont des fibrations principales de fibre un quotient d'un groupe de Lie complexe connexe et abélien.
\end{enumerate}
\end{theo}

Indiquons brièvement comment décrire ces variétés. Le $s^{\textrm{i\`eme}}$ quotient
$$G_s(\ZZ)=\hg{1}{X}/\hg{1}{X}'_s$$
se plonge dans sa complétion de Mal\v{c}ev $G_s(\CC)$ sur laquelle il agit librement et proprement discontinûment\footnote{il faut ici bien se garder de croire que le quotient est compact, ce qui est le cas si on choisit $\RR$ au lieu de $\CC$}. Le quotient $G_s(\ZZ)\backslash G_s(\CC)$ constitue un candidat naturel pour $\albs{s}{X}$ mais il suffit d'examiner la situation pour $s=1$ pour se convaincre qu'il nous faut modifier cette définition naïve. En effet, pour $s=1$, l'abélianisé (sans-torsion) du groupe fondamentale est
$$G_1(\ZZ)=\hg{1}{X}/\hg{1}{X}'_2\simeq \ZZ^{2g}$$
avec $g=q(X)$ et nous avons
$$G_1(\CC)=G_1(\ZZ)\otimes\CC\simeq \CC^{2g}.$$
Le quotient naïf $G_1(\ZZ)\backslash G_1(\CC)$ ne représente pas la variété d'Albanese usuelle et nous devons faire un quotient intermédiaire pour éliminer la "partie anti-holomorphe" de $G_s(\CC)$. Selon R. Hain, les quotients suivants
$$\albs{s}{X}=G_s(\ZZ)\backslash G_s(\CC)/F^0G_s(\CC)$$
vérifient toutes les propriétés énoncées dans le théorème \ref{albanese sup}.

\begin{rem}
Les fibrations $\pi_s:\albs{s+1}{X}\To\albs{s}{X}$ se déduisent aisément des extensions centrales (d'algèbres de Lie et de groupes de Lie respectivement)
\begin{align*}
&0\To \hg{1}{X}_{s+1}/ \hg{1}{X}_{s+2}\otimes\CC\To \mathcal{G}_{s+1}(\CC)\To \mathcal{G}_s(\CC)\To 0\\
\textrm{et }&1\To\hg{1}{X}_{s+1}/ \hg{1}{X}_{s+2}\otimes\CC\To G_{s+1}(\CC)\To G_s(\CC)\To 0
\end{align*}
\end{rem}

\begin{rem}
Outre l'explication fournie ci-dessus, le fait de quotienter $G_s(\CC)$ par le sous-groupe $F^0G_s(\CC)$ a également pour fonction de rendre le morphisme d'Albanese $\alpha_s$ holomorphe.
\end{rem}

Nous aurons également besoin de la propriété suivante des variétés d'Albanese supérieures, conséquence directe de la construction évoquée ci-dessus.
\begin{prop}\label{alb stein}
Pour $s\geq1$, le revêtement universel de $\albs{s}{X}$ est analytiquement isomorphe à un espace $\CC^n$ ; en particulier, $\widetilde{\albs{s}{X}}$ est Stein.
\end{prop}

\section{Convexité holomorphe de $X^{nilp}$}

\subsection{Un cas particulier}\label{cas part}

Nous commençons par exposer un cas particulier intéressant du théorème \ref{shaf nilpotent}, celui correspondant à la situation où l'image de $X$ par son application d'Albanese (usuelle) est lisse. En effet, la démonstration du cas général est assez similaire à celle de ce cas particulier (à la différence qu'il nous faudra prendre en compte toutes les variétés d'Albanese supérieures). Nous allons nous appuyer sur le résultat suivant.

\begin{theo}[\cite{Ca95nilp}]\label{alb lisse}
Soit $X$ une variété kählérienne compacte, $\alpha_X:X\To\alb{X}$ son application d'Albanese et désignons par $Y$ un modèle lisse de l'image $\alpha_X(X)$. Le morphisme naturel $\hg{1}{X}\To\hg{1}{Y}$ induit un morphisme entre les complétions nilpotentes (modulo torsion)
$$\alpha_*^{nilp}:\hg{1}{X}_{nilp}\To\hg{1}{Y}_{nilp}$$
qui est injectif et dont l'image est d'indice finie dans $\hg{1}{Y}_{nilp}$.
\end{theo}
Ce théorème peut également être obtenu en appliquant les résultats de \cite{DGMS} ; la démonstration de \cite{Ca95nilp} procède cependant de façon indépendante et tout à fait élémentaire.\\

\begin{demo}[du théorème \ref{shaf nilpotent} avec $\alpha_X(X)$ lisse]
Soit $Y=\alpha_X(X)\subset\alb{X}$ l'image de $X$ par son application d'Albanese. Le théorème \ref{alb lisse} nous indique que le morphisme $\alpha:X\To Y$ induit un morphisme
$$\alpha_*^{nilp}:\hg{1}{X}_{nilp}\To\hg{1}{Y}_{nilp}$$
injectif et d'image d'indice finie. Considérons alors $\widetilde{Y}$ l'image réciproque de $Y$ dans $\widetilde{\alb{X}}$ :
$$\xymatrix{\widetilde{Y}\ar[d]\ar@{^{(}->}[r] & \widetilde{\alb{X}}\ar[d]\\
Y\ar@{^{(}->}[r] & \alb{X}\,;
}$$
le morphisme $\hg{1}{Y}\To\hg{1}{\alb{X}}$ étant surjectif, $\widetilde{Y}$ est connexe. Comme $\widetilde{Y}$ est fermée dans $\widetilde{\alb{X}}$ qui est Stein, $\widetilde{Y}$ est également une variété de Stein. Le revêtement $\widetilde{Y}\To Y$ est un revêtement abélien ; il est particulier dominé par $Y^{nilp}$ :
$$Y^{nilp}\To \widetilde{Y}\To Y.$$
La variété $Y^{nilp}$ est donc elle aussi Stein, comme revêtement étale d'une variété Stein \cite{St56}. De plus, l'injectivité du morphisme $\alpha_*^{nilp}$ montre que l'application $\alpha:X\To Y$ se relève à $X^{nilp}$ en un diagramme
$$\xymatrix{ X^{nilp}\ar[r]^{\tilde{\alpha}}\ar[d] & Y^{nilp}\ar[d]\\
X\ar[r]^{\alpha} & Y.
}$$
Le fait que $\alpha_*^{nilp}$ ait une image d'indice finie dans $\hg{1}{Y}_{nilp}$ montre immédiatement que l'application $\tilde{\alpha}$ est propre ; $X^{nilp}$ est alors holomorphiquement convexe (car admettant une application propre sur une variété Stein).
\end{demo}

\subsection{Démonstration du théorème \ref{shaf nilpotent}}

Comme il en a été fait mention ci-dessus, la démonstration du cas général procède de la même manière que celle présentée dans le paragraphe précédent. En revanche, il nous faut maintenant prendre en compte toute la tour des variétés d'Albanese et plus seulement $\albs{1}{X}$. Le fait crucial est que cette tour se stabilise à partir d'un certain rang (contrairement à ce qu'il se passe pour les quotients successifs de la suite centrale de $\hg{1}{X}$).

Notons en effet $Y_s=\alpha_s(X)$ les images de $X$ par ses morphismes d'Albanese supérieurs ; nous obtenons de cette sorte un diagramme
\begin{equation*}
\xymatrix{
&&X \ar[ld]_{\alpha_{s+1}}\ar[d]^{\alpha_s}\ar[rrd]^{\alpha_1}&&\\
\ar@{-->}[r]&Y_{s+1}\ar[r]^{\pi_s}&Y_s\ar@{-->}[rr]&&Y_1
}
\end{equation*}
dans lequel les applications $\alpha_s$ et $\pi_s$ sont \emph{surjectives}.
\begin{lem}\label{stabilisation}
Soit $X$ un espace complexe irréductible compact munie d'un diagramme comme ci-dessus ; il existe alors un entier $k\ge1$ à partir duquel les normalisées $\hat{Y_s}$ de $Y_s$ sont biholomorphes entre elles.
\end{lem}

\noindent Tenant pour acquis le lemme précédent, nous pouvons conclure la\\

\begin{demo}[du théorème \ref{shaf nilpotent}]
Soit en effet un entier $k\ge1$ donné par le lemme \ref{stabilisation} et $s\geq k$ ; on sait donc que les normalisées $\hat{Y_s}$ des images $Y_s$ sont toutes biholomorphes à une même variété que nous noterons $\hat{Y}$. On peut alors construire le diagramme suivant
$$\xymatrix{Z_s \ar[r]^{\widetilde{n_s}} \ar[d] & Y'_s\ar@{^{(}->}[r]\ar[d]&\widetilde{\albs{s}{X}}\ar[d] \\
\hat{Y}\ar[r]^{n_s} & Y_s \ar@{^{(}->}[r]&\albs{s}{X}
}$$
dans lequel $Z_s$ est le produit fibré
$$Z_s=\hat{Y}\times_{Y_s} Y'_s.$$
Comme le morphisme $\alpha_s$ se relève en
$$X\stackrel{\hat{\alpha_s}}{\To}\hat{Y}\To\albs{s}{X},$$
on obtient le diagramme suivant
$$\xymatrix{\pi_1(X)\ar@{->>}[rr]^{\alpha_{s*}} \ar[rd]& & \pi_1(\albs{s}{X}) \\
&\pi_1(\hat{Y})\ar[ru]&
}$$
dans lequel la flèche $\alpha_{s*}$ est surjective. On en déduit que le morphisme $\hg{1}{\hat{Y}}\To\hg{1}{\albs{s}{X}}$ est surjectif, ou encore que le produit fibré $Z_s$ est connexe. Le revêtement $Z_s\To \hat{Y}$ est donc un revêtement galoisien (connexe) de groupe de Galois :
\begin{equation}\label{revêtement stein}
\mathrm{Gal}\left(Z_s/\hat{Y}\right)=\hg{1}{\albs{s}{X}}=\hg{1}{X}/\hg{1}{X}'_s.
\end{equation}
Remarquons tout de suite que $Z_s$ est Stein : le morphisme $\widetilde{n_s}$ étant déduit du morphisme de normalisation $n_s$ par produit fibré, il est donc fini et $Y'_s$ est Stein car fermée dans $\widetilde{\albs{s}{X}}$ (proposition \ref{alb stein}). On conclut alors quant au caractère Stein de $Z_s$ car tout revêtement ramifié fini d'un Stein l'est encore \cite[prop. 1.1, p. 252]{FG02}.

On vient de voir (\ref{revêtement stein}) que le groupe $\hg{1}{\hat{Y}}$ admet pour quotient $\hg{1}{X}/\hg{1}{X}'_s$ pour tout $s\ge k$ ; on en déduit donc qu'il admet également pour quotient :
\begin{equation}\label{quotient nilp}
\hg{1}{\hat{Y}}\twoheadrightarrow \hg{1}{X}_{nilp}.
\end{equation}
Notons alors $\widetilde{Y}_\infty\To \hat{Y}$ le revêtement correspondant à (\ref{quotient nilp}). Comme il domine tous les revêtements précédents ($\widetilde{Y}_\infty\To Z_s\To\hat{Y}$), $\widetilde{Y}_\infty$ est encore Stein. Il est alors évident que le morphisme $\hat{\alpha}:X\To \hat{Y}$ se relève en
$$\xymatrix{X^{nilp}\ar[d]\ar[r]^{\tilde{\alpha}} & \widetilde{Y}_\infty\ar[d]\\
X\ar[r]^{\hat{\alpha}}&\hat{Y}
}$$
et que $X^{nilp}$ s'identifie au produit fibré
$$X\times_{\hat{Y}}\widetilde{Y}_\infty.$$
On conclut alors la démonstration puisque $\tilde{\alpha}$ est propre et $\widetilde{Y}_\infty$ est Stein.
\end{demo}

\noindent Il nous reste maintenant à établir la\\

\begin{demo}[du lemme \ref{stabilisation}]
Soit donc un diagramme
\begin{equation}\label{diag}
\xymatrix{
&&X \ar[ld]_{\alpha_{s+1}}\ar[d]^{\alpha_s}\ar[rrd]^{\alpha_1}&&\\
\ar@{-->}[r]&Y_{s+1}\ar[r]^{\pi_s}&Y_s\ar@{-->}[rr]&&Y_1
}
\end{equation}
dans lequel tous les morphismes sont surjectifs. Comme la suite $(\dimm{Y_s})_{s\ge1}$ est croissante et majoré par $\dimm{X}$, on supposera (quitte à renuméroter) que les espaces $Y_s$ ont tous même dimension. 

Nous allons montrer (par récurrence sur la dimension de $X$) que les morphismes $\pi_s$ sont finis à partir d'un certain rang. Si $X$ est une courbe, il n'y a rien à démontrer. Notons alors $d=\dimm{X}-\dimm{Y_1}\ge 0$ et considérons
$$X_1=\set{x\in X\vert\,\dimm{\alpha_1^{-1}(\alpha_1(x))}\geq d+1}.$$
C'est un sous-ensemble analytique fermé \emph{propre} de $X$ ; il est en particulier de dimension strictement inférieure à celle de $X$. Or, si $y\in Y_1\backslash \alpha_1(X_1)$, la commutativité du diagramme :
$$\xymatrix{ & Y_s\ar@{-->}[d]\\
X\ar[ru]^{\alpha_s}\ar[r]^{\alpha_1} & Y_1
}$$
montre que la fibre au dessus de $y$ de l'application composée $Y_s\To Y_1$ doit nécessairement être de dimension nulle (car toute fibre de $X\To Y_s$ est de dimension au moins $d$). On est ramené à restreindre le diagramme (\ref{diag}) à chaque composante irréductible de $X_1$ ; par compacité, il n'y en qu'un nombre fini, toutes de dimensions strictement inférieures à celle de $X$ et on conclut donc grâce à l'hypothèse de récurrence.

La situation finale est donc celle d'un diagramme (\ref{diag}) dans lequel toutes les applications
$$\pi_s:Y_{s+1}\To Y_s$$
sont finies et tous les espaces considérés normaux (quitte à remplacer $X$ et les $Y_s$ par leurs normalisées). Mais, par commutativité du diagramme, le degré de $Y_{s}$ sur $Y_1$ est certainement majoré par une constante indépendante (qui n'est autre que le nombre maximal de composantes connexes des fibres de $X\To Y_1$). Par multiplicativité du degré, on en déduit que le degré de $Y_{s+1}$ sur $Y_s$ doit être égal à 1 à partir d'un certain rang. D'après le théorème principal de Zariski, cela signifie exactement que $\pi_s$ est biholomorphe pour tout $s$ assez grand.
\end{demo}

\vspace*{1cm}

\bibliographystyle{amsalpha}
\bibliography{myref}

\vspace*{0.5cm}
\begin{flushright}
\begin{minipage}{5cm}
Beno\^it CLAUDON\\
Universit\'e Nancy 1\\
Institut Elie Cartan\\
BP 239\\
54 506 Vandoeuvre-l\`es-Nancy\\
Cedex (France)

\vspace*{0.3cm}
\noindent Benoit.Claudon@iecn.u-nancy.fr
\end{minipage}
\end{flushright}

\end{document}